\documentclass[12pt,leqno]{amsart}

\usepackage{amsfonts}
\usepackage{amssymb}
\usepackage{amsmath}
\usepackage{amscd}
\usepackage{amstext}
\usepackage{ifthen}
\usepackage[all]{xy}
\usepackage{enumerate}
\usepackage{fancyhdr}   
\usepackage[pdftex]{color}

\makeatletter
\def\@cite#1#2{{\m@th\upshape\bfseries%
[{#1\if@tempswa{\m@th\upshape\mdseries, #2}\fi}]}}
\makeatother

\newtheorem{theorem*}{Theorem}
\newtheorem{theorem}{Theorem}[section]
\newtheorem{lemma}[theorem]{Lemma}
\newtheorem{corollary}[theorem]{Corollary}
\newtheorem{proposition}[theorem]{Proposition}

\theoremstyle{definition}

\newtheorem{remark}[theorem]{Remark}

\newtheorem{conjecture}{Conjecture}

\numberwithin{equation}{section}

  \newcommand{\A}{{\mathcal{A}}}
  \newcommand{\B}{{\mathcal{B}}}

  \newcommand{\E}{{\mathcal{E}}}
  \newcommand{\F}{{\mathcal{F}}}

  \newcommand{\M}{{\mathcal{M}}}
  \newcommand{\N}{{\mathcal{N}}}

  \newcommand{\X}{{\mathcal{X}}}
  \newcommand{\Y}{{\mathcal{Y}}}
  \newcommand{\Z}{{\mathcal{Z}}}

\renewcommand{\phi}{\varphi}
\newcommand{\upchi}{{\raise.35ex\hbox{\ensuremath{\chi}}}}

\newcommand{\ad}{\operatorname{Ad}}

\newcommand{\di}{\operatorname{d}}

\newcommand{\id}{{\operatorname{id}}}

\newcommand{\m}{\operatorname{m}}
\newcommand{\n}{\operatorname{n}}

\newcommand{\li}{\operatorname{-lim}}

\begin{document}

\title[]{Dual operator algebras close to injective von Neumann algebras}

\author[J. Roydor]{Jean Roydor}

\thanks{The author is supported by JSPS}

\begin{abstract}
We prove that if a non-selfadjoint dual operator algebra admitting a normal virtual diagonal and an injective von Neumann algebra are close enough for the Kadison-Kastler's metric, then they are similar. The bound explicitly depends on the norm of the normal virtual diagonal. This is inspired from E. Christensen's work on perturbation of operator algebras.
\end{abstract}

\date{}
\maketitle

\section{Introduction}

Perturbation theory of operator algebras in the sense of Kadison-Kastler is an active area at the moment with two recent remarkable articles \cite{C5} and \cite{C6}.
The starting point of this paper is the conjunction of perturbation theory of operator algebras and a conjecture on amenable non-selfadjoint operator algebras. Let us first recall this conjecture and a dual version of it, then we will explain the connection with our main result.\\
\indent A conjecture (raised by G. Pisier) asserts that a non-selfadjoint amenable operator algebra $\A$ should be similar to a nuclear $C^*$-algebra i.e. there is an invertible operator $S$ such that $S\A S^{-1}$ is a $C^*$-algebra. In his memoir \cite{J0}, B.E. Johnson characterized amenability for Banach algebras by the existence of a virtual diagonal. Recall that injectivity for von Neumann algebras can be characterized by the existence of a normal virtual diagonal (in the sense of E.G. Effros \cite{E}, see Subsection \ref{Pre:d} below for more details). Therefore, a dual version of the preceding conjecture would be:
\begin{conjecture}A unital dual operator algebra $\M$ admitting a normal virtual diagonal $u$ should be similar to an injective von Neumann algebra.
 \end{conjecture}
 In that case, it is expected that the similarity constant is controlled by a non-decreasing function of the norm of the normal virtual diagonal. Both of these conjectures are still open (see the introduction of \cite{M} for more details).\\
\indent In 1972, R.V. Kadison and D. Kastler defined a metric $d$ on the collection of all subspaces of the bounded operators on a fixed Hilbert space (see Subsection \ref{Pre:p}). They conjectured that sufficiently close $C^*$-algebras are necessarily unitarily conjugated (see \cite{KK}). A great amount of work around this conjecture has been achieved since then (see \cite{C5} for a nice introduction on this topic). Especially by E. Christensen, he proved it for the class of type I von Neumann algebras (in \cite{C0}) and for the class of injective von Neumann algebras (in \cite{C1}, \cite{C2}). Very recently, Kadison-Kastler's conjecture has been proved for the class of separable nuclear $C^*$-algebras in \cite{C5} by E. Christensen, A.M. Sinclair, R.R. Smith, S.A. White and W. Winter (see also \cite{C4}). Let us state E. Christensen's first result on perturbation of injective von Neumann algebras (this result has subsequently been improved in \cite{C2}):
\begin{theorem*}(E. Christensen, Th. 4.1 \cite{C1})\label{T:C} Let $\M, \N$ be two von Neumann subalgebra of a fixed $\mathbb{B}(H)$. We suppose that $\M$ has Schwartz's property $(P)$ and $\N$ has the extension property.\\ If $d(\M,\N)<1/169$, then there is a unitary $U$ in the von Neumann algebra generated by $\M \cup \N$ such that $U \M U^*=\N$. Moreover, $\Vert U-I_H \Vert \leq 19 d(\M,\N)^{1/2}$.
\end{theorem*}
  After the work of A. Connes \cite{Co1}, \cite{Co2} and U. Haagerup \cite{H}, we know that Schwartz's property $(P)$, the extension property and injectivity (and thus the existence of a normal virtual diagonal) are equivalent conditions for von Neumann algebras.\\
 \indent Conjecture 1 mentioned above leads to the following question: can we replace, in the preceding theorem, $\M$ by a unital non-selfadjoint dual operator algebra admitting a normal virtual diagonal ? In other words, is the selfadjointness hypothesis on $\M$ necessary ? Indeed, assume for a moment that Conjecture 1 is true, then there would be an invertible $S$ such that $S\M S^{-1}$ is an injective von Neumann algebra (moreover, $$d(\M,S\M S^{-1}) \leq 2(1+\Vert S \Vert \Vert S^{-1} \Vert)\Vert S - I_H \Vert$$ and this last quantity is controlled by a non-decreasing function of $\Vert u \Vert$).
If $d(\M,\N)$ is small enough such that the following strict inequality holds $$d(\N,S\M S^{-1}) \leq d(\N,\M) + d(\M,S\M S^{-1}) < \frac{1}{169},$$ then (from Theorem \ref{T:C} above) the injective von Neumann algebras $\N$ and $S\M S^{-1}$ would be unitarily conjugated, so $\M$ and $\N$ would be similar. Therefore, it is not totally incongruous to try to replace $\M$ by a unital dual operator algebra admitting a normal virtual diagonal.\\
\indent In this paper, we prove (see Theorem \ref{T:jr5}):

\begin{theorem*}\label{T:jr3} Let $\M,\N \subset \mathbb{B}(H)$ be two unital $w^*$-closed operator algebras. Suppose that $\M$ admits a normal virtual diagonal $u$ and $\N$ is an injective von Neumann algebra.\\
If $d(\M,\N) < \frac{1}{656\Vert u \Vert}$, then there exists an invertible operator $S$ in the $w^*$-closed algebra generated by $\M \cup \N$ such that $S\M S^{-1}=\N$. Moreover, $\Vert S-I_H \Vert \leq 656\Vert u \Vert d(\M,\N)$.
\end{theorem*}

 Note that von Neumann algebras enjoy a self-improvement phenomenon, if a von Neumann algebra admits a normal virtual diagonal then it admits a normal virtual diagonal of norm one, see \cite{H}, \cite{E} or \cite{EK} (self-improvement phenomena are frequent for selfadjoint algebras, for instance nuclearity constant and exactness constant). This may explain why in Theorem \ref{T:C} the bound is a universal constant, whereas in Theorem \ref{T:jr3}, the bound depends on the feature of the non-selfadjoint algebra involved. Moreover, from Theorem 7.4.18 (1) in \cite{BLM} and Remark \ref{R:1} below, if a unital dual operator algebra admits a normal virtual diagonal of norm one, then it is necessarily a von Neumann algebra (no similarity is needed in this extreme case). Hence, Theorem \ref{T:C} corresponds exactly to the case $\Vert u \Vert$ equals 1 in Theorem \ref{T:jr3} (as the unitary $U$ is obtained by taking the polar decomposition of $S$, see Lemma 2.7 in \cite{C0}). Our bound in this special case is not as good as E. Christensen's one, but the important point is that we have removed the selfadjointness hypothesis on $\M$. This is not a minor modification, knowing that non-selfadjoint algebras are less rigid than selfadjoint ones (no order structure for instance) and fewer tools are available (no continuous or Borel functional calculus), so our proof requires new ingredients as we shall see immediately.\\
 \indent Let us sketch the main lines of our proof, there are three steps (as in E. Christensen's work \cite{C1}):
\begin{enumerate}[Step 1]
\item find a linear isomorphism, between the two algebras, which is close to the identity representation,
\item find an algebra homomorphism which is close to the previous linear isomorphism,
\item prove that this algebra homomorphism is similar to the identity representation.
\end{enumerate}
For the first step, as $\N$ is injective, one just has to take the restriction to $\M$ of a completely contractive projection onto $\N$. This gives a linear isomorphism $T$ from $\M$ onto $\N$ which is close to the identity representation of $\M$. But actually, there is an extra effort between the first and second step: one is working with dual operator algebras, so the linear isomorphism must be $w^*$-continuous (in order to apply certain averaging procedure afterwards). For this, E. Christensen used Tomiyama's decomposition into normal/singular part of bounded linear maps defined on von Neumann algebras. But when $\M$ is non-selfadjoint, such decomposition is not available. Hence, we have to consider the normal part of $T^{-1}$. This $w^*$-continuous linear isomorphism from $\N$ onto $\M$ is not necessarily completely positive, moreover the target algebra $\M$ is not necessarily self-adjoint, thus we can not use E. Christensen's averaging trick (see Lemma 3.3 in \cite{C1}) to accomplish the second step. The idea is to turn to Banach algebras results. More precisely, we will use a dual operator space version of a B.E. Johnson theorem on almost multiplicative maps (see \cite{J}). Indeed, the issue here is that we need to preserve the $w^*$-continuity, but we can not use the normal projective tensor product of dual Banach spaces (as we could not check its associativity, see Section \ref{Bej}). This second step will force us to work with the normal Haagerup tensor product of dual operator spaces. Finally, the third step which is related to a more general problem on neighboring representations (already mentioned in \cite{KK}), is done by an averaging technique. However, because of the second step, we have had to work in the operator space category and as a consequence we had to assume that the algebras $\mathbb{M}_n(\M)$ nearly embed in $\mathbb{M}_n(\N)$ uniformly in $n$ (see the notion of near cb-inclusion defined in Subsection \ref{Pre:p} below). As an intermediate result, we prove a perturbation theorem with a near cb-inclusion assumption (see Theorem \ref{T:jr4}). Therefore, our final task is to notice that the existence of a normal virtual diagonal is an ``automatic near cb-inclusion" condition (see Lemma \ref{L:3}).

\section{Preliminaries}\label{Pre}

For background on completely bounded maps, operator space theory and non-selfadjoint algebra theory, the reader is referred to \cite{BLM}, \cite{ER1}, \cite{Pa} and \cite{P}. Especially Section 2.7 in \cite{BLM} for background on dual operator algebras.

\subsection{Perturbation theory}\label{Pre:p}
We first recall definitions and notations commonly used in perturbation theory of operator algebras (see e.g. \cite{C6}).
Let $H$ be an Hilbert space, and $\mathbb{B}(H)$ be the von Neumann algebra of all bounded operators on $H$. Let $\E,\F$ be two subspaces of $\mathbb{B}(H)$. We denote by $d$ the Kadison-Kastler's metric i.e. $d(\E,\F)$ denotes the Hausdorff distance between the unit balls of $\E$ and $\F$. More explicitly,
\begin{align*}d(\E,\F)=\inf \{ \gamma > 0 :& \forall x \in B_\E, \exists x' \in B_\F, \Vert x - x' \Vert < \gamma ~\mbox{and}~\\
                      &        \forall y \in B_\F, \exists y' \in B_\E, \Vert y - y' \Vert < \gamma                  \},
\end{align*}
where $B_\E$ (respectively $B_\F$) denotes the unit ball of $\E$ (respectively $\F$).\\
 Let $\gamma>0$, we write $\E \subseteq^\gamma \F$  if for any $x$ in the unit ball of $\E$, there exists $y$ in $\F$ such that $$\Vert x-y \Vert \leq \gamma.$$ We also write $\E \subset^\gamma \F$ if there exists $\gamma'<\gamma$ such that $\E \subseteq^{\gamma'} \F$.\\ We will also need the notion of \textit{near cb-inclusion}. As usual in operator space theory, $\mathbb{M}_n(\E)$, the subspace of $n \times n$ matrices with coefficients in $\E$ is normed by the identification $\mathbb{M}_n(\E) \subset \mathbb{M}_n(\mathbb{B}(H))=\mathbb{B}(\ell^2_n \otimes H)$. We write
$$\E \subseteq_{cb}^\gamma \F$$ if, for all $n$, $\mathbb{M}_n(\E) \subseteq^\gamma  \mathbb{M}_n(\F)$.

\subsection{The normal projective tensor product and the normal Haagerup tensor product}\label{Pre:t}
  For $\X$, $\Y$ dual operator spaces, we denote $(\X \otimes_h \Y)_\sigma^*$ the space of all completely bounded bilinear forms which are separately $w^*$-continuous (see Paragraph 1.5.4 of \cite{BLM} for definition of completely bounded bilinear maps). The normal Haagerup tensor product, denoted $\otimes_{\sigma h}$, can be defined as \begin{equation}\label{eqp1}\X \otimes_{\sigma h} \Y= ((\X \otimes_h \Y)_\sigma^*)^*,\end{equation} see Paragraph 1.6.8 of \cite{BLM}. The normal Haagerup tensor product is characterized by the following universal property: $\X \otimes \Y$ is $w^*$-dense in $\X \otimes_{\sigma h} \Y$ and for any dual operator space $\Z$, for any $w^*$-continuous completely contractive bilinear map $B:\X \times \Y \to \Z$, there exists a (unique) $w^*$-continuous completely contractive linear map $\tilde{B}:\X \otimes_{\sigma h} \Y \to \Z$ such that $\tilde{B}(x \otimes y)=B(x,y)$, for all $x \in \X$, $y \in \Y$.\\
We will also need the normal projective tensor product $\widehat{\otimes}_\sigma$ of dual Banach spaces. If $\X$, $\Y$ are dual Banach spaces, $$\X \widehat{\otimes}_{\sigma} \Y= ((\X \widehat{\otimes} \Y)_\sigma^*)^*,$$
where $(\X \widehat{\otimes} \Y)_\sigma^*$ denotes the space of all bounded bilinear forms on $\X \times \Y$ which are separately $w^*$-continuous. The normal projective tensor product enjoys similar universal property as the normal Haagerup tensor product, but for separately $w^*$-continuous bounded bilinear maps instead of separately $w^*$-continuous completely bounded (for von Neumann algebras, the projective normal tensor product appeared for instance in \cite{E} under the name binormal projective tensor product).\\
These two tensor products are ``functorial'' in the sense that, if $L_i:\X_i \to \Y_i$ ($i=1,2$) are bounded (resp. completely bounded) $w^*$-continuous linear maps between dual Banach spaces (resp. dual operator spaces), then there is a unique bounded (resp. completely bounded) $w^*$-continuous linear map $$L_1 \widehat{\otimes}_\sigma L_2:\X_1 \widehat{\otimes}_\sigma \X_2 \to \Y_1 \widehat{\otimes}_\sigma \Y_2$$ (resp. $L_1 \otimes_{\sigma h} L_2:\X_1 \otimes_{\sigma h} \X_2 \to \Y_1 \otimes_{\sigma h} \Y_2$) extending $L_1 \otimes L_2$. Moreover, $\Vert L_1 \widehat{\otimes}_\sigma L_2 \Vert \leq \Vert L_1 \Vert \Vert L_2 \Vert$ (resp. $\Vert L_1 \otimes_{\sigma h} L_2 \Vert_{cb} \leq \Vert L_1 \Vert_{cb} \Vert L_2 \Vert_{cb}$).\\
The main difference between these two tensor products is that the normal Haagerup tensor product is associative (see Lemma 2.2 in \cite{BK}), whereas the normal projective tensor product does not seem to be associative in general (this difference will have important consequence for us in Section \ref{Bej}).

\subsection{Normal virtual diagonals and normal virtual $h$-diagonals}\label{Pre:d} Normal virtual diagonals appeared implicitly in \cite{H} and explicitly in \cite{E} (see p. 147). In this paper, we also need the notion normal virtual $h$-diagonal (called reduced normal virtual diagonal in \cite{E}, see also Paragraph 7.4.8 of \cite{BLM} for more details). Let us just recall this notion. Replacing the normal Haagerup tensor product by the normal projective tensor product in the following, one can obtain analogously the definition of normal virtual diagonal.\\
Let $\M$ be a unital dual operator algebra, let us recall the $\M$-bimodule structure of $\M \otimes_{\sigma h} \M$. Let
$\psi \in (\M \otimes_{ h} \M)_\sigma^*$ and $a,b,c,d \in \M$ $$\langle b \cdot \psi \cdot a, c \otimes d \rangle=\psi(ac,db).$$ Hence by duality, one can define actions of $\M$ on $\M \otimes_{\sigma h} \M= ((\M \otimes_h \M)_\sigma^*)^*$. One can check that these actions are determined on the elementary tensors by: $$a \cdot (c \otimes d)\cdot b= ac \otimes db.$$
On a dual operator algebra, the multiplication is a separately $w^*$-continuous completely contractive bilinear map (see (1) Proposition 2.7.4 in \cite{BLM}). Consequently, it induces a $w^*$-continuous complete contraction $$\m_{\sigma h}:\M \otimes_{\sigma h} \M \to \M.$$
A \textit{normal virtual $h$-diagonal} for $\M$ is an element $u \in \M \otimes_{\sigma h} \M$ satisfying:
\begin{enumerate}[(C1)]
\item for any $m \in \M$, $m \cdot u= u \cdot m$,
\item $\m_{\sigma h}(u)=1$.
\end{enumerate}
 Note that condition (C2) implies that the norm of a normal virtual $h$-diagonal is always greater or equal to 1.\\
\begin{remark}\label{R:1}
 Note that the inclusion $(\X \otimes_h \Y)_\sigma^* \subset  (\X \widehat{\otimes} \Y)_\sigma^*$ induces, by duality, a contraction from $\M \widehat{\otimes}_\sigma \M$ into $\M \otimes_{\sigma h} \M$ and this contraction sends normal virtual diagonals into normal virtual $h$-diagonals. Consequently, if $\M$ admits a normal virtual diagonal, it admits a normal virtual $h$-diagonal.
\end{remark}


\section{A B.E. Johnson's theorem revisited}\label{Bej}

The aim of this section is to find a solution to the second step mentioned in the Introduction section. In \cite{J}, B.E. Johnson proved that an approximately multiplicative map defined on an amenable Banach algebra is close to an actual algebra homomorphism. His result is the Banach algebraic version of an earlier result due to D. Kazhdan for amenable groups (see \cite{K}).
If $L$ is a linear map between operator algebras $\M$ and $\N$, we denote $L^{\vee}:\M \times \M \to \N$ the bilinear map defined by $$L^{\vee}(x,y)=L(xy)-L(x)L(y).$$ This enables us to measure the defect of multiplicativity of $L$.\\
In our present case, we have to take into account the dual operator space structure of our algebras, starting from a $w^*$-continuous linear map from $\M$ into $\N$, we have to obtain $w^*$-continuous algebra homomorphism. This will force us to work in the category of operator spaces. The proof of Theorem 3.1 in \cite{J} is by induction, the algebra homomorphism is the limit (in operator norm) of a sequence of linear maps with multiplicativity defect tending to zero. The problem is that these linear maps are defined using the $w^*$-topology of the target algebra (see equation $(*)$ in the proof of Theorem 3.1 in \cite{J}). Here, to justify the $w^*$-continuity of these linear maps, we have to consider a trilinear map defined on $\M \times \M \times \M$ (see Equation (\ref{b0eq}) below). But the normal projective tensor product does not seem associative. To circumvent this difficulty, we will instead work with normal Haagerup tensor product, which is associative (see Lemma 2.2 in \cite{BK}). As a consequence, we have to control the cb-norm of the bilinear map $L^\vee$.

\begin{remark}\label{R:0} Actually, this difficulty concerning the associativity of the normal projective tensor product has already been encountered in disguise. The main issue in \cite{E} is that one can not check whether the Banach $\M$-bimodule $\M \widehat{\otimes}_\sigma \M$ is normal or not. But if one assumes that the normal projective tensor product is associative, then it is not too difficult to see that $\M \widehat{\otimes}_\sigma \M$ is normal.
\end{remark}

\begin{theorem}\label{T:j} Let $\M,\N$ be two unital dual operator algebras. We suppose that $\M$ has a normal virtual $h$-diagonal $u \in \M \otimes_{\sigma h} \M$. Then, for any $\varepsilon \in ]0,1[$, for any $\mu >0$, there exists $\delta >0$ such that: for every unital $w^*$-continuous linear map $L:\M \to \N$ satisfying $\Vert L \Vert_{cb} \leq \mu$ and $\Vert L^{\vee} \Vert_{cb} \leq \delta$, there is a unital $w^*$-continuous completely bounded algebra homomorphism $\pi:\M \to \N$ such that $\Vert L - \pi \Vert_{cb} \leq \varepsilon$.
\end{theorem}
\begin{proof}
Let $\varepsilon \in ]0,1[, \mu >0$ and $L$ be a unital $w^*$-continuous linear map from $\M$ into $\N$ such that $\Vert L \Vert_{cb} \leq \mu$.\\
The trilinear map \begin{equation}\label{b0eq}(x,y,z) \in \M \times \M \times \M \mapsto L(x)L^\vee(y,z) \in \N\end{equation}
is separately $w^*$-continuous and completely bounded. By the universal property of the normal Haagerup tensor product, it extends to a $w^*$-continuous completely bounded linear map $\Lambda_L:\M \otimes_{\sigma h} \M \otimes_{\sigma h} \M \to \N$ such that $$\Lambda_L(x \otimes y \otimes z)= L(x)L^\vee(y,z)$$ and $\Vert \Lambda_L \Vert_{cb} \leq \Vert L  \Vert_{cb} \Vert L^\vee  \Vert_{cb}$. By definition and associativity of the normal Haagerup tensor product, the linear map
$$m \in \M \longmapsto u \otimes m \in (\M \otimes_{\sigma h} \M) \otimes_{\sigma h} \M=\M \otimes_{\sigma h} \M \otimes_{\sigma h} \M$$ is $w^*$-continuous (see Equation (\ref{eqp1})). We can define $R:\M \to \N$ by $$R(m)=\Lambda_L(u \otimes m)$$ which is $w^*$-continuous and \begin{equation}\label{b1eq}\Vert R \Vert_{cb} \leq \Vert u \Vert \Vert L  \Vert_{cb} \Vert L^\vee  \Vert_{cb}.\end{equation}
As $u \in \M \otimes_{\sigma h} \M$, there is a net $(u_t)_t$ in $\M \otimes \M$ converging to $u$ in the $w^*$-topology of $\M \otimes_{\sigma h} \M$. For any $t$, there are finite families $(a_k^t)_k,(b_k^t)_k$ of elements in $\M$ such that $$u_t=\sum_k a_k^t \otimes b_k^t.$$ Now fix $m \in \M$, once again by definition and associativity of the normal Haagerup tensor product, the linear map
$$v \in \M \otimes_{\sigma h} \M \longmapsto v \otimes m \in (\M \otimes_{\sigma h} \M) \otimes_{\sigma h} \M=\M \otimes_{\sigma h} \M \otimes_{\sigma h} \M$$ is $w^*$-continuous as well. Hence, using the $w^*$-continuity of $\Lambda_L$, we obtain $$R(m)=w^*\li _t \Lambda_L(u_t \otimes m)=w^*\li _t \sum_k L(a_k^t)L^\vee(b_k^t,m).$$
From this point, we just need to check that computations of \cite{J} Theorem 3.1 remain valid with matrix coefficients.
 Fix $n \in \mathbb{N}$, let $x,y$ in the unit ball of $\mathbb{M}_n(\M)$ (in the following computation, $I_n$ denotes the identity matrix in $\mathbb{M}_n$, and the others subscripts $n$ denote the $n^{th}$ ampliation of a linear or bilinear map), then as in \cite{J}, we have:

\begin{align*}
 (L+R)_{n}^\vee (x,y)&=\\
 &L_{n}^\vee (x,y) - w^*\li _t \sum_k L_{n}(I_n \otimes a_k^tb_k^t)L_{n}^\vee(x,y)\\
 &-R_{n}(x)R_{n}(y)\\
 &+ w^*\li _t \sum_k L_{n}^\vee(x,I_n \otimes a_k^t)L_{n}^\vee(I_n \otimes b_k^t,y)\\
 &+ w^*\li _t \sum_k L_{n}^\vee(I_n \otimes a_k^t,I_n \otimes b_k^t)L_{n}^\vee(x,y)\\
 &+ w^*\li _t \sum_k L_{n}(I_n \otimes a_k^t)L_{n}((I_n \otimes b_k^t)xy)-L_{n}(x(I_n \otimes a_k^t))L_{n}((I_n \otimes b_k^t)y)\\
 &- w^*\li _t \sum_k (L_{n}(I_n \otimes a_k^t)L_{n}((I_n \otimes b_k^t)x)-L_{n}(x(I_n \otimes a_k^t))L_{n}(I_n \otimes b_k^t))L_{n}(y)
\end{align*}
To evaluate the norm of $(L+R)_{n}^\vee$, we treat each line of the right hand-side successively.\\
As $u$ is a normal virtual $h$-diagonal, $w^*\li _t \sum_k a_k^tb_k^t=1$. But $L$ is unital and $w^*$-continuous, so $w^*\li _t \sum_k L_{n}(I_n \otimes a_k^tb_k^t)=1$ and the first line of the right hand-side is $0$. Clearly, the norm of the term in the second line is bounded by $\Vert R \Vert_{cb}^2$. Considering the bilinear map $$(a,b) \in \M \times \M \mapsto L_{n}^\vee(x,I_n \otimes a)L_{n}^\vee(I_n \otimes b,y) \in \mathbb{M}_n(\N),$$ we can see that
that the norm of the term in the third line is bounded by $\Vert u \Vert \Vert L^\vee \Vert_{cb}^2$. Similarly, the norm of term in the fourth line is bounded by $\Vert u \Vert \Vert L^\vee \Vert_{cb}^2$. For the term in the fifth line, note that its $(i,j)$-entry is
$$w^*\li _t \sum_k \sum_{p=1}^n (L(a_k^t)L(b_k^tx_{ip}y_{pj})-L(x_{ip}a_k^t)L(b_k^ty_{pj})) \in \N.$$
But $u$ is a normal virtual $h$-diagonal, so for any $i,p$, $$w^*\li _t (x_{ip} \cdot u_t-u_t \cdot x_{ip})=0,$$
hence for any $i,j,p$,
 $$w^*\li _t (\sum_k x_{ip} \cdot a_k^t \otimes b_k^t \cdot y_{pj}- \sum_k a_k^t \otimes b_k^t \cdot x_{ip}y_{pj})=0.$$
 The bilinear map $$(x,y) \in \M \times \M \mapsto L(x)L(y) \in \N$$ extends to a $w^*$-continuous map, consequently the term in the fifth line is $0$. Analogously, the term in the sixth line is $0$ as well. Finally we obtain that \begin{equation}\label{b2eq}\Vert (L+R)^\vee \Vert_{cb} \leq 2\Vert u \Vert \Vert L^\vee \Vert_{cb}^2+ \Vert R \Vert_{cb}^2 \leq (2\Vert u \Vert +\Vert u \Vert^2 \Vert L  \Vert_{cb}^2) \Vert L^\vee  \Vert_{cb}^2 \end{equation}
 Now we are in position to follow the induction of \cite{J} with cb-norms instead of norms (for the reader convenience, we reproduce it here), the important point is that each $L^q$ (thus each $R^q$) defined below are $w^*$-continuous. Define \begin{equation}\label{b3eq}\delta=\frac{\varepsilon}{4 \Vert u \Vert + 8\mu^2\Vert u \Vert^2}.\end{equation} Suppose that $\Vert L^\vee \Vert_{cb} \leq \delta$.
 Inductively, we define a sequence a linear maps from $\M$ into $\N$ by $L_0=L$, $R_0=R$ and for $q \geq 0$, $$L^{q+1}=L^q+R^q \quad \mbox{and} \quad R^{q+1}(\cdot)=\Lambda_{L^{q+1}}(u \otimes \cdot).$$ We also define $\mu_q=(2-2^{-q})\mu$, $\delta_q=2^{-q}\delta$. By induction, we prove that for all $q$, $\Vert (L^{q})^\vee \Vert_{cb} \leq \delta_q$ and $\Vert L^{q} \Vert_{cb} \leq \mu_q$. It's obvious for $q=0$. Then using the inequality (\ref{b2eq}) above, we have $$\Vert (L^{q+1})^\vee \Vert_{cb} \leq (2\Vert u \Vert +\Vert u \Vert^2 \mu_q^2) \delta_q^2 \leq \delta_{q+1}$$
 and using (\ref{b1eq}) to majorize the cb-norm of $R^q$, we obtain
 $$ \Vert L^{q+1} \Vert_{cb} \leq \mu_q + \Vert u \Vert \mu_q \delta_q \leq \mu_{q+1},$$
 (the last inequality coming from the fact that $\Vert u \Vert \delta \leq 4^{-1}$).
 Consequently, $$\Vert R^q \Vert_{cb} \leq \Vert u \Vert \Vert L^q  \Vert_{cb} \Vert (L^q)^\vee  \Vert_{cb} \leq 2\Vert u \Vert \mu \delta_q,$$ so $\sum_{q \geq 0} R^q$ converges in cb-norm. We can define $$\pi=L+\sum_{q \geq 0} R^q,$$ in other words $\pi=\lim_q L^q$, so $\pi$ is $w^*$-continuous. Hence $\pi^\vee=\lim_q (L^q)^\vee$, but we proved that $\Vert (L^q)^\vee \Vert_{cb} \leq \delta_q$, so $\pi$ is multiplicative. Moreover, $\Vert \pi-L \Vert_{cb}= \Vert \sum_{q \geq 0} R^q \Vert_{cb} \leq 4 \Vert u \Vert \mu \delta < \varepsilon.$
\end{proof}

\begin{remark} One important point which does not appear in the statement of the previous theorem is that $\delta$ is an explicit function of $\mu$, $\varepsilon$ and $\Vert u \Vert$ (see Equation (\ref{b3eq}) above).
\end{remark}

\section{Neighboring representations}

We now show that two representations of a dual operator algebra, admitting a normal virtual $h$-diagonal, which are close enough in cb-norm are necessarily similar. Apparently, this phenomena is well-known to Banach algebraists (see e.g. \cite{R0} Chapter 8). We give here a quick proof for dual operator algebras. This proposition will enable us to perform the third step mentioned in the Introduction section.\\
If $S \in \mathbb{B}(H)$ is an invertible operator, we denote $\ad_S$ the similarity implemented by $S$.

\begin{proposition}\label{P:rep} Let $\M$ be a unital dual operator algebra. We suppose that $\M$ has a normal virtual $h$-diagonal $u \in \M \otimes_{\sigma h} \M$. Let $\pi_1, \pi_2$ be two unital $w^*$-continuous completely bounded representations on the same Hilbert space $K$.\\
If $\Vert \pi_1 - \pi_2 \Vert_{cb}<\Vert u \Vert ^{-1}  \max \{ \Vert \pi_1  \Vert_{cb} ^{-1}, \Vert \pi_2  \Vert_{cb} ^{-1} \}$, then there exists an invertible operator $S$ in the $w^*$-closed algebra generated by $\pi_1(\M) \cup \pi_2(\M)$ such that $\pi_1=\ad_S \circ \pi_2$.
\end{proposition}
\begin{proof} Let $\pi_1,\pi_2$ as above. For two completely bounded $w^*$-continuous linear maps $F,G:\M \to \mathbb{B}(K)$, we denote (with notation of Section \ref{Pre:t}) $$\Psi_{F,G}=\m_{\sigma h} \circ (F \otimes_{\sigma h} G),$$ which is a completely bounded $w^*$-continuous linear map defined on $\M \otimes_{\sigma h} \M$.
 Now, define $$S=\Psi_{\pi_1,\pi_2}(u) \in \mathbb{B}(K).$$ As $u \in \M \otimes_{\sigma h} \M$, there is a net $(u_t)_t$ in $\M \otimes \M$ converging to $u$ in the $w^*$-topology of $\M \otimes_{\sigma h} \M$. For any $t$, there are finite families $(a_k^t)_k,(b_k^t)_k$ of elements in $\M$ such that $$u_t=\sum_k a_k^t \otimes b_k^t.$$ Hence, $S=w^*\li _t \sum_k \pi_1(a_k^t)\pi_2(b_k^t).$
  Let $m \in \M$, then
   \begin{align*}\pi_{1}(m)S=& \pi_{1}(m).w^*\li _t \sum_k \pi_1(a_k^t)\pi_2(b_k^t)\\
=&w^*\li _t \sum_k \pi_1(ma_k^t)\pi_2(b_k^t) \\
 = & w^*\li _t \Psi_{\pi_1,\pi_2}(m\cdot u_t)\\
 = & \Psi_{\pi_1,\pi_2}(m\cdot u)
 \end{align*}
 Analogously, we can show that $$S \pi_2(m)=\Psi_{\pi_1,\pi_2}(u \cdot m).$$
 But $u$ is a normal virtual $h$-diagonal, so $m\cdot u=u \cdot m$, hence $$\pi_{1}(m)S=S \pi_2(m).$$

To conclude, we just need to prove that $S$ is invertible. Without loss of generality we can assume that $\Vert \pi_1 - \pi_2 \Vert_{cb}<\Vert u \Vert ^{-1} \Vert \pi_1  \Vert_{cb} ^{-1}$. As above, we have $\Psi_{\pi_1,\pi_1}(u)= w^*\li _t \sum_k \pi_1(a_k^t)\pi_1(b_k^t)$. Using the condition (C2) defining a normal virtual $h$-diagonal, we obtain,
\begin{align*}\Psi_{\pi_1,\pi_1}(u)=& w^*\li _t \sum_k \pi_1(a_k^tb_k^t)\\
=&\pi_1(w^*\li _t \sum_k a_k^tb_k^t) \\
 = & \pi_1(\m_{\sigma h}(u))\\
 =& \pi_1(1)\\
 =&I_K.\end{align*}
Consequently,

\begin{align*}\Vert S - I_K \Vert=& \Vert \Psi_{\pi_1,\pi_2}(u)- \Psi_{\pi_1,\pi_1}(u) \Vert\\
=&\Vert \Psi_{\pi_1,\pi_2-\pi_1}(u) \Vert \\
 \leq & \Vert u \Vert \Vert \pi_1 - \pi_2 \Vert_{cb} \Vert \pi_1  \Vert_{cb}\\
  < & 1\end{align*}
 which ends the proof.
\end{proof}

\section{Proof of the main Theorems}

We start this section with a very simple lemma that we will use repeatedly in the proof of the next theorem, we just sketch the proof. Recall that
$T^{\vee}$ denotes the bilinear map from $\M \times \M$ into $\N$ defined by $T^{\vee}(x,y)=T(xy)-T(x)T(y)$. Also in this section, we denote $\id_\A$ the identity representation of a concretely represented operator algebra $\A$.

\begin{lemma}\label{L:0}
Let $\A,\B \subset \mathbb{B}(H)$ be two operator algebras and $T:\A \to \B$ be a completely bounded linear map. Then,
\begin{enumerate}[(i)]
\item $\Vert T^{\vee} \Vert_{cb} \leq (2+\Vert T \Vert_{cb})\Vert T - \id_\A \Vert_{cb}$.
\item If $\Vert T - \id_\A \Vert_{cb} < 1$, then $T$ is injective and has closed range. Moreover, if there exists $\alpha \in [0,1[$ such that for any $y$ in the unit ball of $\B$, there is $x$ in $\A$ satisfying $\Vert T(x)-y \Vert \leq \alpha$, then $T$ is bijective and $$\Vert T^{-1} \Vert_{cb} \leq \frac{1}{1-\Vert T - \id_\A \Vert_{cb}}.$$
   \end{enumerate}
     \end{lemma}
     \begin{proof}
     Let $x,y$ in the unit ball of $\mathbb{M}_n(\M)$, then $(i)$ follows from the decomposition
\begin{align*}
(T^\vee)_n(x,y)&=T_n(xy)-T_n(x)T_n(y)\\
&=T_n(xy)-xy+xy-xT_n(y)+ xT_n(y)- T_n(x)T_n(y).
\end{align*}
  The first assertion of $(ii)$ follows from \begin{align*}\Vert T_n(x) \Vert &\geq \vert \Vert T_n(x) - x \Vert-\Vert x \Vert \vert\\
&\geq (1-\Vert T - \id_\A \Vert_{cb})\Vert x \Vert.
\end{align*}
The surjectivity of $T$ is proved by induction. Let $y$ in the unit ball of $\N$, for any integer $j$, we can find $t_1,\dots,t_j$ in the range of $T$ such that
$$ \Vert y-(t_1+t_2+ \cdots +t_j) \Vert \leq \alpha^j.$$
As $\alpha < 1$, we conclude that $y$ belongs to the closure of the range of $T$.
     \end{proof}

Note that, in the following theorem, $\M$ is just assumed to be a dual operator algebra, but we require a near cb-inclusion of $\M$ into $\N$ (see Subsection 2.1).

\begin{theorem}\label{T:jr4} Let $\M,\N \subset \mathbb{B}(H)$ be two unital $w^*$-closed operator algebras. We suppose that $\N$ is injective von Neumann algebra.\\
If $\N \subset^1 \M$ and $\M \subseteq_{cb}^{\gamma} \N$, with $\gamma < 1/164$, then there exists an invertible operator $S$ in the $w^*$-closed algebra generated by $\M \cup \N$ such that $S\M S^{-1}=\N$.
\end{theorem}
\begin{proof}
Since $\N$ is injective, there is a completely contractive projection $P$ from $\mathbb{B}(H)$ onto $\N$. Denote $T=P_{\vert \M}$. Let $x$ in the unit ball of $\mathbb{M}_n(\M)$, then there is $y$ in $\mathbb{M}_n(\N)$ such that $\Vert x - y \Vert \leq \gamma$.
$$\Vert T_n(x) - x \Vert =\Vert T_n(x-y) - (x-y) \Vert_{cb} \leq 2\gamma,$$ hence $$\Vert T - \id_\M \Vert_{cb} \leq 2\gamma<1.$$ Let's prove that $T$ is surjective. Since $\N \subset^1 \M$, there is $\gamma' < 1$ such that $\N \subseteq^{\gamma'} \M$. Let $y$ in the unit ball of $\N$, then there exists $x$ in $\M$ such that $\Vert y-x \Vert \leq \gamma'$, so $$\Vert T(x)-y \Vert= \Vert T(x)-T(y) \Vert \leq \gamma',$$ therefore from Lemma \ref{L:0}$(ii)$, $T$ is a linear cb-isomorphism and \begin{equation}\label{e0}\Vert T^{-1} \Vert_{cb} \leq \frac{1}{1-2\gamma}.\end{equation}
The problem is that $T$ is not necessarily $w^*$-continuous, so we are going to consider the normal of $T^{-1}$ (see \cite{T}, we denote with an exponent $\n$ the normal part of a linear map defined on $\N$). Note first that \begin{equation}\Vert T^{-1} - \id_\N \Vert_{cb} \leq \Vert T^{-1} \Vert_{cb} \Vert T - \id_\M \Vert_{cb} \leq \frac{2\gamma}{ 1-2\gamma}.\end{equation}
Let $V=(T^{-1})^{\n}:\N \to \M$ the normal part of $T^{-1}$ . Using Lemma \ref{L:0}$(ii)$ again, let's show that $V$ is a completely bounded $w^*$-continuous linear isomorphism from $\N$ onto $\M$. As taking the normal part is a completely contractive operation, we have
\begin{equation}\label{e2}\Vert V-\id_\N \Vert_{cb} = \Vert (T^{-1}-\id_\N)^{\n} \Vert_{cb} \leq \Vert T^{-1}-\id_\N \Vert_{cb} \leq \frac{2\gamma}{ 1-2\gamma},\end{equation}
 thus $V$ is injective and has closed range. Let $y$ in the unit ball of $\M$. Pick $x$ in $\N$ such that $\Vert x-y \Vert \leq \gamma$. Thus $\Vert x \Vert \leq 1 + \gamma$ and
\begin{align*}
\Vert V(x)-y \Vert & \leq \Vert V(x)-x \Vert+\Vert x-y \Vert\\
& \leq \frac{2\gamma}{1-2\gamma}(1+\gamma)  + \gamma \\
& \leq \frac{5\gamma}{1-2\gamma}<1
\end{align*}
so $V$ is surjective.\\
In order to apply Theorem \ref{T:j}, we need to unitize $V$. From Equation (\ref{e2}), $\Vert V(1)-1 \Vert \leq \frac{2\gamma}{ 1-2\gamma}<1$, $V(1)$ is invertible in $\M$ and we obtain
\begin{equation}\label{e3}\Vert V(1)^{-1} \Vert \leq \frac{1}{1-\Vert V(1)-1 \Vert} \leq \frac{1-2\gamma}{ 1-4\gamma}.\end{equation} Denote $L=V(1)^{-1}V$, so $L$ is a unital $w^*$-continuous completely bounded isomorphism from $\N$ onto $\M$ and from Equation (\ref{e0}) we have \begin{equation}\label{e4}\Vert L \Vert_{cb} \leq \Vert V(1)^{-1} \Vert \Vert V \Vert_{cb} \leq \Vert V(1)^{-1} \Vert \Vert T^{-1} \Vert_{cb} \leq \frac{1}{ 1-4\gamma}.\end{equation}
Let's compute the defect of multiplicativity of $L$.
Fix $n$, let $x$ in unit ball of $\mathbb{M}_n(\N)$, from Equations (\ref{e2}) and (\ref{e3}) we obtain
\begin{align*}
\Vert L_n(x)-x \Vert & \leq \Vert I_n \otimes V(1)^{-1}(V_n(x)-x) \Vert+\Vert I_n \otimes V(1)^{-1}x-x \Vert\\
& \leq \Vert V(1)^{-1} \Vert \Vert V-\id_\N \Vert_{cb}+ \Vert V(1)^{-1} \Vert\Vert V(1)-1 \Vert\\
& \leq \frac{4\gamma}{1-4\gamma}
\end{align*}
which means that \begin{equation}\label{e5}\Vert L - \id_\N \Vert_{cb}\leq \frac{4\gamma}{1-4\gamma}.\end{equation}
Therefore, by Lemma \ref{L:0}$(i)$ and Equation (\ref{e4}) we obtain
$$\Vert L^{\vee} \Vert_{cb} \leq (2+\Vert L \Vert_{cb})\Vert L - \id_\N \Vert_{cb}\leq \frac{12\gamma}{(1-4\gamma)^2}.$$
We want to apply Theorem \ref{T:j} to $L$.  Put $$\mu=\frac{1}{ 1-4\gamma}~~\mbox{and}~~ \delta=\frac{12\gamma}{(1-4\gamma)^2}.$$
As $\N$ is an injective von Neumann algebra, we can find a normal virtual $h$-diagonal $u$ of norm $1$ (see \cite{E}, \cite{EK}), thus (see Equation (\ref{b3eq})) let $$\varepsilon=\delta(4 \Vert u \Vert + 8\mu^2\Vert u \Vert^2)=\frac{12\gamma}{(1-4\gamma)^2}(4+\frac{8}{(1-4\gamma)^2}).$$
    We can then apply Theorem \ref{T:j} to $L$ and find a unital $w^*$-continuous completely bounded homomorphism $\pi:\N \to \M$ such that $$\Vert L - \pi \Vert_{cb} \leq \varepsilon.$$
 Consequently, from Equation (\ref{e5}) \begin{align*}
\Vert \pi-\id_\N \Vert_{cb} & \leq \Vert \pi-L \Vert_{cb}+\Vert L-\id_\N \Vert_{cb}\\
& \leq \varepsilon + \frac{4\gamma}{1-4\gamma}
\end{align*}
and this last quantity is strictly smaller than 1, because $\gamma<1/164$. Therefore, we can apply Proposition \ref{P:rep} to $\pi$ and $\id_\N$ and find an invertible operator $S$ in the $w^*$-closed algebra generated by $\M \cup \N$ such that $$\ad_S \circ \pi=\id_\N$$ (in particular $\pi$ is injective and has closed range). To achieve the proof, it is sufficient to prove that the range of $\pi$ is $\M$. Let $y$ is in the unit ball of $\M$, then $$\Vert \pi(L^{-1}(y))-y \Vert \leq \Vert \pi-L \Vert_{cb} \Vert L^{-1} \Vert_{cb},$$ so by Lemma \ref{L:0}$(ii)$, we just need to check that this last quantity is strictly smaller than one. From Equation (\ref{e5}) $$\Vert L^{-1} \Vert_{cb} \leq \frac{1}{1-\Vert L - \id_\N \Vert_{cb}} \leq \frac{1-4\gamma}{1-8\gamma},$$
it follows that $$\Vert \pi-L \Vert_{cb} \Vert L^{-1} \Vert_{cb} \leq \frac{1-4\gamma}{1-8\gamma}\varepsilon$$
which is strictly smaller than 1, because $\gamma<1/164$.
\end{proof}

At this point, we want to get rid of the near cb-inclusion hypothesis appearing in the previous theorem. The question is to find conditions of ``automatic near cb-inclusion" on the algebra $\M$. More explicitly, under which condition a near inclusion $\M \subseteq^\gamma \N$ implies automatically a near cb-inclusion ? For $C^*$-algebras, E. Christensen isolated property $D_k$ which ensures such ``automatic near cb-inclusion" result (see \cite{C2}). Recall that a $C^*$-algebra $\A$ has property $D_k$ if for any unital $*$-representation $(\pi,K)$ one has: $$\forall x \in \mathbb{B}(K),\quad \di(x,\pi(\A)') \leq k \Vert \delta (x)_{\vert \pi(\A)} \Vert,$$
where $\di$ denotes the usual distance between between subsets and $\delta (x)$ denotes the inner derivation implemented by $x$ on $\mathbb{B}(K)$. It is well-known that amenable $C^*$-algebras (or injective von Neumann algebras) have $D_1$, the next easy lemma is the non-selfadjoint analog of this fact (it also works for amenable Banach algebras).

\begin{lemma}\label{L:2} Let $\M$ be a unital dual operator algebra admitting a normal virtual diagonal $u \in \M \widehat{\otimes}_{\sigma} \M$.\\ Then, for any unital $w^*$-continuous contractive representation $(\pi,K)$ of $\M$ satisfying $\pi(\M)=\overline{\pi(\M)}^{w^*}$, we have: \begin{equation}\label{e6}\forall x \in \mathbb{B}(K),\quad \di(x,\pi(\M)') \leq \Vert u \Vert \Vert \delta (x)_{\vert \pi(\M)} \Vert.\end{equation}
\end{lemma}
\begin{proof}
Let's denote $\N=\pi(\M) \subset \mathbb{B}(K)$ and $v=\pi \widehat{\otimes}_{\sigma} \pi (u) \in \N \widehat{\otimes}_{\sigma} \N$, hence $\Vert v \Vert \leq \Vert u \Vert$. Since $\pi$ has $w^*$-closed range, $v$ is a normal virtual diagonal for the dual operator algebra $\N$. Note that $\mathbb{B}(K)$ is obviously a normal dual Banach $\N$-bimodule (in the sense of \cite{R0} Definition 4.4.6). Now, let $x \in \mathbb{B}(K)$ and consider the $w^*$-continuous bounded derivation $D=\delta(x)_{\vert \N}: \N \to \mathbb{B}(K)$. Adapting the proof of B.E Johnson's Theorem on characterization of amenability by virtual diagonals, we know that there is $\phi \in \mathbb{B}(K)$ such that $D=\delta(\phi)_{\vert \N}$. Actually $\phi=D \widehat{\otimes}_{\sigma} \id_\N (v)$, so $\Vert \phi \Vert \leq \Vert D \Vert \Vert v \Vert$. As $D=\delta(x)_{\vert \N}=\delta(\phi)_{\vert \N}$, $x-\phi \in \N'$. Therefore,
$$\di(x,\pi(\M)')=\di(\phi,\N') \leq \Vert \phi \Vert \leq \Vert D \Vert \Vert v \Vert \leq \Vert u \Vert \Vert \delta (x)_{\vert \pi(\M)} \Vert,$$ which ends the proof
\end{proof}

\begin{lemma}\label{L:3} Let $\M \subset \mathbb{B}(H)$ be a unital $w^*$-closed operator algebra admitting a normal virtual diagonal $u \in \M \widehat{\otimes}_{\sigma} \M$. Let $\N$ be an injective von Neumann subalgebra of $\mathbb{B}(H)$.\\ Then, for any $\gamma>0$, $\M \subseteq^\gamma \N$ implies $\M \subseteq_{cb}^{4 \Vert u \Vert \gamma} \N$.
\end{lemma}
\begin{proof}
This follows from the previous lemma and the first lines of the proof of Theorem 3.1 in \cite{C2}, with $D=\mathbb{M}_n$ (for arbitrary $n$), with $k=\Vert u \Vert$ (by Equation (\ref{e6}) and the $3/2$ must replaced by $1$ because $\N$ is injective, so we obtain $4 \Vert u \Vert \gamma$ instead of $6 k \gamma$.
\end{proof}

Using the previous Lemma and Theorem \ref{T:jr4} above, we can deduce the main result:

\begin{theorem}\label{T:jr5} Let $\M,\N \subset \mathbb{B}(H)$ be two unital $w^*$-closed operator algebras. Suppose that $\M$ admits a normal virtual diagonal $u \in \M \widehat{\otimes}_{\sigma} \M$. We suppose that $\N$ is an injective von Neumann algebra.\\
If $\N \subset^1 \M$ and $\M \subseteq^{\gamma} \N$, with $\gamma < \frac{1}{656\Vert u \Vert}$, then there exists an invertible operator $S$ in the $w^*$-closed algebra generated by $\M \cup \N$ such that $S\M S^{-1}=\N$. Moreover, $\Vert S-I_H \Vert \leq 656\Vert u \Vert \gamma$.
\end{theorem}

This question of ``automatic near cb-inclusion" can be thought as an analog of the ``automatic complete boundedness" question for homomorphisms (or equivalently Kadison's similarity problem). For this problem, G. Pisier defined the notion of length for operator algebras (see \cite{P1}, \cite{P2}). The connection between this notion of length and property $D_k$ is now well-known for $C^*$-algebras, see \cite{C6}. As we are working with dual operator algebras, C. Le Merdy's notion of length (or degree) denoted $d_*$ in \cite{LM} is more appropriate (we call this quantity \textit{normal length} in the following corollary).

\begin{corollary} Let $\M,\N \subset \mathbb{B}(H)$ be two unital $w^*$-closed operator algebras. Suppose that $\M$ has finite normal length at most $d_*$ with constant at most $C>0$. We suppose that $\N$ is an injective von Neumann algebra.\\
If $\N \subset^1 \M$ and $\M \subseteq^{\gamma} \N$, with $\gamma < (1+\frac{1}{164C})^{1/d_*}-1$, then there exists an invertible operator $S$ in the $w^*$-closed algebra generated by $\M \cup \N$ such that $S\M S^{-1}=\N$ and consequently, $d_*(\M)\leq 2$.
\end{corollary}
\begin{proof} If $\M \subseteq^\gamma \N$, then $\M \subseteq_{cb}^{C((1+\gamma)^{d_*}-1)} \N$ as in Proposition 2.10 in \cite{C6}. The consequence follows from the similarity degree characterization of injectivity for von Neumann algebra in \cite{P3}.
\end{proof}

\begin{remark}
One can notice that, in our proof, operator space theory is needed only for Step 2, this is also the case in E. Christensen's work (Lemma 3.3 \cite{C1} or Theorem 4.1 in \cite{C2} where Stinespring factorization is used). It seems that operator space theory can not be avoided in those perturbation questions. Note also that our new ingredient for this Step 2 (an operator space version of B.E. Johnson's theorem on approximately multiplicative maps) has cohomological interpretation, see \cite{J}. In section 3, we switched from the normal projectif tensor product to the normal Haagerup tensor product, this actually corresponds to a change from the bounded cohomology to the completely bounded cohomology. Then we used an ``automatic complete boundness" argument for near inclusions. The same strategy is used in \cite{C3.5}.
\end{remark}

\textbf{Acknowledgements.} The author would like to thank Stuart White for introducing him to perturbation theory of operator algebras and suggesting him to look at non-selfadjoint algebras. The author is grateful to Narutaka Ozawa for hosting him at the University of Tokyo.

\email{Jean Roydor, Department of Mathematical Sciences, Tokyo, 153-8914, Japan.\\
roydor@ms.u-tokyo.ac.jp}


\begin{thebibliography}{99}

\bibitem{BK} D.P. Blecher, U. Kashyap, \textit{Morita equivalence of dual operator algebras}. J. Pure Appl. Algebra  212  (2008),  no. 11, 2401--2412.
\bibitem{BLM} D.P. Blecher, C. Le Merdy, \textit{Operator algebras and their modules---an operator space approach}. London Mathematical Society Monographs. New Series, 30. Oxford Science Publications.

       \bibitem{C0} E. Christensen, \textit{Perturbations of type I von Neumann algebras}. J. London Math. Soc. (2)  9  (1974/75), 395--405.
\bibitem{C1} E. Christensen, \textit{Perturbations of operator algebras}. Invent. Math. 43  (1977), no. 1, 1--13.
\bibitem{C2} E. Christensen, \textit{Near inclusions of $C^{\ast}$-algebras}. Acta Math. 144 (1980), no. 3-4, 249--265.
\bibitem{C3} E. Christensen, \textit{Extensions of derivations. II}. Math. Scand.  50  (1982), no. 1, 111--122.
\bibitem{C3.5} E. Christensen, F. Pop, A. M. Sinclair, R.R. Smith \textit{Hochschild Cohomology of Factors with Property Gamma}. Ann. Math.
158 (2003), 635--659.
\bibitem{C4} E. Christensen, A.M. Sinclair, R.R. Smith, S.A. White, W. Winter, \textit{The spatial isomorphism problem for close separable nuclear $C^*$-algebras}. Proc. Natl. Acad. Sci. USA 107 (2010), no. 2, 587--591.
\bibitem{C5} E. Christensen, A.M. Sinclair, R.R. Smith, S.A. White, W. Winter, \textit{Perturbations of nuclear $C^*$-algebras}. Acta Mathematica, to appear. Preprint Arxiv math.OA/0910.4953v1.
\bibitem{C6} E. Christensen, A.M. Sinclair, R.R. Smith, S.A. White, W. Winter, \textit{Perturbations of C*-algebraic invariants}. Geom. Funct. Anal., 20 (2010) 368-397.
\bibitem{Co1} A. Connes, \textit{Classification of injective factors. Cases $II_{1},$ $II_{\infty },$ $III_{\lambda },$ $\lambda \not=1$}.  Ann. of Math. (2)  104  (1976), no. 1, 73115.
\bibitem{Co2} A. Connes, \textit{On the cohomology of operator algebras}. J. Functional Analysis 28 (1978), no. 2, 248253.
\bibitem{EP} G.K. Eleftherakis, V.I. Paulsen, \textit{Stably isomorphic dual operator algebras}. Math. Ann. 341 (2008), no. 1, 99--112.
\bibitem{E} E.G. Effros, \textit{Amenability and virtual diagonals for von Neumann algebras}. J. Funct. Anal. 78 (1988), no. 1, 137--153.
\bibitem{EK} E.G. Effros, A. Kishimoto, \textit{Module maps and Hochschild-Johnson cohomology}. Indiana Univ. Math. J. 36 (1987), no. 2, 257--276.
\bibitem{ER1} E.G. Effros, Z.J. Ruan \textit{Operator spaces}. London Mathematical Society Monographs. New Series, 23. The Clarendon Press, Oxford University Press, New York, 2000.
        \bibitem{H} U. Haagerup, \textit{A new proof of the equivalence of injectivity and hyperfiniteness for factors on a separable Hilbert space}.  J. Funct. Anal.  62  (1985),  no. 2, 160201.
    \bibitem{J0} B.E. Johnson, \textit{Cohomology in Banach algebras}. Mem. Amer. Math. Soc. 127 (1972), 1-96.
\bibitem{J} B.E. Johnson, \textit{Approximately multiplicative maps between Banach algebras}. J. London Math. Soc. (2)  37  (1988),  no. 2, 294--316.
\bibitem{KK} R.V. Kadison, D. Kastler \textit{Perturbations of von Neumann algebras. I. Stability of type}.  Amer. J. Math.  94  (1972), 3854.
\bibitem{K} D. Kazhdan, \textit{On $\varepsilon $-representations}.  Israel J. Math.  43  (1982), no. 4, 315--323.
\bibitem{LM} C. Le Merdy, \textit{The weak$^*$ similarity property on dual operator algebras}.  Integral Equations Operator Theory  37  (2000),  no. 1, 7294.
\bibitem{M} L.W. Marcoux, \textit{triangularizable, total reduction algebras}.  J. Lond. Math. Soc. (2)  77  (2008),  no. 1, 164182.
\bibitem{Pa} V.I. Paulsen, Vern \textit{Completely bounded maps and operator algebras}. Cambridge Studies in Advanced Mathematics, 78. Cambridge University Press, Cambridge, 2002.
\bibitem{PS1} V.I. Paulsen, R.R. Smith, \textit{Multilinear maps and tensor norms on operator systems}. J. Funct. Anal.  73  (1987),  no. 2, 258--276.
\bibitem{P1} G. Pisier, \textit{The similarity degree of an operator algebra}. St. Petersburg Math. J.  10  (1999),  no. 1, 103146.
\bibitem{P2} G. Pisier, \textit{The similarity degree of an operator algebra. II}.  Math. Z.  234  (2000),  no. 1, 5381.
\bibitem{P} G. Pisier, \textit{Introduction to operator space theory}. London Mathematical Society Lecture Note Series, 294. Cambridge University Press, Cambridge, 2003.
\bibitem{P3} G. Pisier, \textit{A similarity degree characterization of nuclear $C^*$-algebras}. Publ. Res. Inst. Math. Sci. 42 (2006), no. 3, 691704.
\bibitem{R0} V. Runde, \textit{Lectures on amenability}. Lecture Notes in Mathematics, 1774. Springer-Verlag, Berlin, 2002.
    \bibitem{T} J. Tomiyama, \textit{On the projection of norm one in $W\sp*$-algebras. III}. T\^ohoku Math. J. (2)  11  1959 125129.
\end{thebibliography}
\end{document}